\documentclass[12pt]{article}
\usepackage{enumerate}
\usepackage{amssymb,amsmath,amsfonts,amsthm,mathrsfs}
\usepackage[colorlinks=true, pdfstartview=FitV, linkcolor=blue, citecolor=blue, urlcolor=blue]{hyperref}
\usepackage{graphicx,color}
\usepackage{cite}
\usepackage{indentfirst}

\allowdisplaybreaks

\begin{document}
\def\rn{{\mathbb R^n}}  \def\sn{{\mathbb S^{n-1}}}
\def\co{{\mathcal C_\Omega}}
\def\z{{\mathbb Z}}
\def\nm{{\mathbb (\rn)^m}}
\def\mm{{\mathbb (\rn)^{m+1}}}
\def\n{{\mathbb N}}
\def\cc{{\mathbb C}}

\newtheorem{defn}{Definition}
\newtheorem{thm}{Theorem}
\newtheorem{lem}{Lemma}
\newtheorem{cor}{Corollary}
\newtheorem{rem}{Remark}

\title{\bf\Large Sharp constants for a class of linear integral operators on higher-dimensional product spaces
\footnotetext{
\newline\indent\hspace{1mm} {\it 2020 Mathematics Subject Classification}: Primary 42B25; Secondary 42B20, 47H60, 47B47.}}

\date{}
\author{Xiang Li, Zunwei Fu\footnote{Corresponding author} and Zhongci Hang\footnote{Corresponding author}}
\maketitle
\begin{center}
\begin{minipage}{13cm}
{\small {\bf Abstract:}\quad
In this paper, we will study a class of linear integral operators with the nonnegative kernels on higher-dimensional product spaces, the norms of the operators can be obtained by integral of the product of the kernel function and finitely many basic functions. As application, we obtain the sharp constant for  linear Hilbert operator on Higher-dimensional product spaces.
 }
\end{minipage}
\end{center}

\section[Introduction]{Introduction}
How to compute the sharp constants of some inequalities is a very important question in harmonic analysis. For a class of linear integral operators with the nonnegative kernels,  which satisfy some homogeneity and the rotational invariance condition,Beckner \cite{Beckner2} obtained the following conclusion which is also denoted by Stein-Weiss Lemma.

\begin{lem}\label{main_1}
	Suppose that $K$ is a nonnegative kernel defined on $\mathbb{R}^{n}\times \mathbb{R}^{n}$ continuous on any domain that
	excludes the point $(0, 0)$, homogeneous of degree $-n$.\,
	\begin{equation}\label{10}
		K(\delta u, \delta v) = \delta ^{-n}K(u, v)
	\end{equation}
	and
	\begin{equation}\label{11}
		K(Ru, Rv) = K(u, v)
	\end{equation}
	for any $ R \in \mathrm SO(n)$. Then $K$ defines an integral operator
	$$
	T f(x)=\int_{\mathbb{R}^n} K(x, y) f(y) d y
	$$
	which maps $L^p(\mathbb{R}^n)$ to $L^p(\mathbb{R}^n)$ $i.e.$,
	\begin{equation}
		\|T f\|_{L^p(\mathbb{R}^n)} \leq C\|f\|_{L^p(\mathbb{R}^n)}
	\end{equation}
	holds for $1<p<\infty$, where the optimal constant is given by
	\begin{equation}
		C=\int_{\mathbb{R}^n} K(e_1, y)|y|^{-\frac{n}{p^{\prime}}} d y,
	\end{equation}
	where $e_1=(1,0,\cdots,0)$ is a unit vector in the first coordinate direction in $\mathbb{R}^n$, and $\mathrm SO(n)$ denotes a set of the rotation transformation on $\mathbb{R}^n$.
\end{lem}
The inequality (\ref{main_1}) was given by Stein and Weiss \cite{Stein}, then Beckner  \cite{Beckner1, Beckner2, Beckner3} pointed out that the constant $C$  is the best constant and gave the sharp constant for the Hardy-Littlewood-Sobolev inequality with some power weights £¬more information readers can refer to \cite{Benyi, DW} . Our main result is inspired by Lemma \ref{main_1}. For more information about higher-dimensional product spaces, the reader can refer to \cite{Shimo, He, Lusz}. Now , we give the definition of a class of linear integral operator with nonnegative kernel on higher-dimensional product spaces.

\begin{defn}
	Let $m \in \mathbb{N}, n_i \in \mathbb{N}, x_i \in \mathbb{R}^{n_i}, 1 \leq i \leq m$, and $f$ be a nonnegative measurable function on $\mathbb{R}^{n_1} \times \mathbb{R}^{n_2} \times \cdots \times \mathbb{R}^{n_m}$.Then operator $T_m $ with the kernel $K$ is defined by
	\begin{equation}
		T_mf(x)=\int_{\mathbb{R}^{n_1}}\cdots\int_{\mathbb{R}^{n_m}} K(x_1, y_1)\cdots K(x_m, y_m) f(y_1, y_2\cdots y_m) d y_1 d y_2\cdots d y_m.
	\end{equation}
	
\end{defn}
It is not difficult to find that our work is actually to generalize the Stein-Weiss Lemma in higher-dimensional product space. Now, we give our main  result.
\begin{thm}\label{main_2}
	Let $1<p<\infty$, $m \in \mathbb{N}, n_i \in \mathbb{N}, x_i \in \mathbb{R}^{n_i}, 1 \leq i \leq m$. If $f \in L^p (\mathbb{R}^{n_1}\times\mathbb{R}^{n_2}\times\cdots\times\mathbb{R}^{n_m} )$, then we have
	\begin{equation}
		\|T_m f\|_{L^p(\mathbb{R}^{n_1} \times \mathbb{R}^{n_2} \times \cdots \times \mathbb{R}^{n_m})} \leq C\|f\|_{L^p(\mathbb{R}^{n_1} \times \mathbb{R}^{n_2} \times \cdots \times \mathbb{R}^{n_m})}.
	\end{equation}
	Moreover,
	\begin{equation}
		\|T_m f\|_{L^p(\mathbb{R}^{n_1} \times \mathbb{R}^{n_2} \times \cdots \times \mathbb{R}^{n_m})\rightarrow L^p(\mathbb{R}^{n_1} \times \mathbb{R}^{n_2} \times \cdots \times \mathbb{R}^{n_m})} = C,
	\end{equation}
	where
	\begin{equation}
		C=\int_{\mathbb{R}^{n_1}}\cdots\int_{\mathbb{R}^{n_m}} K(e_1, y_1)\cdots K(e_m, y_m) |y_1|^{-\frac{n_1}{p}} |y_2|^{-\frac{n_2}{p}} \cdots |y_m|^{-\frac{n_m}{p}} d y_1 d y_2\cdots d y_m.
	\end{equation}
\end{thm}
To prove above Theorem, we first provide some definitions and Lemmas which will be used in the following.Some Lemmas are well known in many classic books so we omit their proofs.

\section{Auxiliary lemmas}
\begin{defn}
	A topological group $\mathcal{G}$ is  Hausdorff topological spaces that is also a group with law $(x, y) \rightarrow x y$ such that the maps $(x, y) \rightarrow x y$ and $x \rightarrow x^{-1}$ are continuous.
\end{defn}
\begin{defn}
	A topological group $\mathcal{G}$ is locally compact, if each point has a compact neighborhood.
\end{defn}
\begin{lem}\label{main_3}
	A topological group $\mathcal{G}$ is locally compact if and only if the identity element has a compact
	neighborhood.
\end{lem}
\begin{lem}\label{main_4}
	The product space $X = X_1\times \cdots \times X_n $ is a locally compact topological group if and only
	if each $X_i, i = 1, 2, \ldots, n$ ,is a locally compact topological group.
\end{lem}
\begin{lem}\label{main_5}
	Suppose that $\mathcal{G}$ is a locally compact topological group. Then $\mathcal{G}$ possesses a positive measure $\mu$ on the Borel sets such that $\mu$ is nonzero on all nonempty open sets and is left invariant, which means $\mu(s A)=\mu(A)$ for all $s \in \mathcal{G}$ and  any measurable set A. Furthermore, left Haar measure is unique up to positive multiplicative constants. A locally compact topological group $\mathcal{G}$ with a Haar measure $\mu$ is denoted by $(\mathcal{G}, \mu)$.
\end{lem}
The proofs of Lemmas \ref{main_3}-\ref{main_5} can be found in \cite{Gra, Lang} .

\begin{lem}\label{main_6}
	(1) The multiplicative group $\mathbb{R}^{+}=(0, \infty)$ is a locally compact group, and $G=(0, \infty)^n$.\\
	(2) For the multiplicative group $G$ with group operation
	$$
	x y=(x_1 y_1, \ldots, x_n y_n),
	$$
	the Haar measure $\mu$ on $G$ is $\frac{d x}{x_1 x_2 \cdots x_n}$, and $\mu$ is a $\sigma$-finite measure.\\
	(3) If $A^{-1}={x^{-1}: x \in A}$, then $\mu(A^{-})=\mu(A)$.
\end{lem}
\begin{proof}
	The proof of the first assertion of Lemma \ref{main_6} is trivial.We omit it here.
	Let us turn to the proof of $(2)$ and $(3)$. For all $c=(c_1, c_2, \ldots, c_n) \in G$, set $y=(y_1, y_2, \ldots, y_n)$ and $y_i=$ $c_i x_i, i=1,2, \ldots, n$. We have
	$$
	d y=\prod_{i=1}^n c_i d x_i.
	$$
	For any measurable set $A$ in $G$, we have
	$$
	\mu(c A)=\int_{c A} \frac{d x}{x_1 \cdots x_n}=\int_A \frac{d y}{y_1 \cdots y_n}=\mu(A),
	$$
	where we use the change of variables $x=c y$.
	Assume that $A_m=\left(\frac{1}{m}, m\right)^n, m \in \mathbb{Z}^{+}$, then $\mu(A_m)<\infty$. Since
	$$
	G=\bigcup_{m=1}^{\infty}\left(\frac{1}{m}, m\right)^n ,
	$$
	then $\mu$ is a $\sigma$-finite measure on G. This has finished the proof of (2).
	
	Assume that $0 \leq c_i<d_i \leq \infty, i=1,2, \ldots, n, 1 / 0=\infty$, and $1 / \infty=0$. If $A=(c_1, d_1) \times \cdots \times(c_n, d_n)$, then $A^{-}=(d_1^{-}, c_1^{-}) \times \cdots \times(d_n^{-}, c_n^{-})$. By using the change of variables $x=c y$, we have
	$$
	\mu(A^{-})=\int_{A^{-}} \frac{d x}{x_1 \cdots x_n}=\int_A \frac{d y}{y_1 \cdots y_n}=\mu(A) .
	$$
	By the classical method in measure theory, one concludes $\mu(A^{-})=\mu(A)$  holds for any Borel measurable $A$ in $G$. This has finished the proof of (3).
\end{proof}
\begin{defn}
	Define
	$$
	L^p(\mathcal{G}, \mu)=\left\{f:\|f\|_{L^p(\mathcal{G}, \mu)}=\left(\int_{\mathcal{G}}|f(x)|^p d \mu(x)\right)^{1 / p}<+\infty\right\}, \quad 1 \leq p<+\infty,
	$$
	when $p=\infty$, we have
	$$
	L^{\infty}(\mathcal{G}, \mu)={f:\|f\|_{L^{\infty}}<+\infty}.
	$$
	It is well known that H\"older's inequality holds in $L^p(\mathcal{G}, \mu)$. Moreover, Fubini's theorem also holds, provided that $\mu$ is a $\sigma$-finite measure on $\mathcal{G}$. In this paper, $L^p(G, \mu)$ denotes the function space on the locally compact group $G=(0, \infty)^n$ with respect to the Haar measure $\mu=\frac{d x}{x_1 x_2 \cdots x_n}$. Others such as $L^p(\mathbb{R}^{+}), L^p(\mathbb{R}^n)$ and $L^p(G)$ are ones with respect to the Lebesgue measure.
\end{defn}
\begin{lem}\label{8}(see\cite{Gra})
	Let $(X,\mu)$ and $ (Y,v)$ be two $\sigma$-finite measure spaces.Let $T$ be a linear operator defined  the set of all finitely simple functions on X and taking values in the set of measurable functions on $Y$.Let $1\leq p_0, p_1, q_0, q_1 \leq \infty$ and  we assume that
	$$
	\begin{aligned}
		& \|T(f)\|_{L^{q_0}} \leq M_0\|f\|_{L^{p_0}}, \\
		& \|T(f)\|_{L^{q_1}} \leq M_1\|f\|_{L^{p_1}}
	\end{aligned}
	$$
	hold for all finitely simple functions $f$ on $X$. Then for all $0<\theta<1$ , we have
	$$
	\|T(f)\|_{L^q} \leq M_0^{1-\theta} M_1^\theta\|f\|_{L^p}
	$$
	hold for all finitely simple functions $f$ on $X$, where
	\begin{equation}\label{main_7}
		\frac{1}{p}=\frac{1-\theta}{p_0}+\frac{\theta}{p_1} \quad \text { and } \quad \frac{1}{q}=\frac{1-\theta}{q_0}+\frac{\theta}{q_1} \text {. }
	\end{equation}
	Consequently, when $p<\infty$, by density, $T$ has a unique bounded extension from $L^p(X, \mu)$ to $L^q(Y, v)$ where $p$ and $q$ are as in Theroem \ref{main_7}.
\end{lem}
\begin{lem}{\label{main_13}}
	Let $(X,\mu)$ and$(Y,v)$ be two $\sigma$-finite measure spaces. Let $T$ be a linear operator defined on the set of all finitely simple functions on $X$ and taking values in the set of measurable functions on $Y$. If $T$ satisfies
	$$
	\|T\|_{L^{\infty}(X, \mu)\rightarrow L^{\infty}(Y, v)} \leq C
	$$
	and
	$$
	\|T\|_{L^{1}(X, \mu)\rightarrow L^{1}(Y, v)} \leq C,
	$$
	then the operator $T$ maps $L^{p}$ to $L^p$, and
	$$
	\|T\|_{L^{p }\rightarrow L^{p}} \leq C,
	$$
	for $1\leq p \leq\infty$ with
	$$
	0\leq \frac{1}{p}\leq 1 .
	$$
\end{lem}
\begin{lem}
	Suppose that a function $K: \mathbb{R}^{n} \rightarrow[0, \infty)$, and each $f$ is a measure function on $\mathbb{R}^n$. Define linear integral operator as
	$$
	\widetilde{T}_K (f)(x)=: \int_{\mathbb{R}^{n}} K(y-x) f(y) d y,
	$$
	where $x,  y \in \mathbb{R}^n$. If $1 \leq p \leq \infty$  and
	$$
	0 \leq \frac{1}{p} \leq 1,
	$$
	then we have
	\begin{equation}\label{9}
		\|\widetilde{T}_K\|_{L^{p}  \rightarrow L^p}=\|K\|_{L^1(\mathbb{R}^{n })}.
	\end{equation}
\end{lem}
\begin{cor}\label{main_19}
	Suppose that a function $H: \mathbb{G}^{n } \rightarrow[0, \infty)$, and  $f$ is a measure function on $\mathbb{G}^n$ with respect to Haar measure. Define linear integral operator as
	$$
	\widetilde{T}_H(f)(x):=\int_{\mathbb{G}^{n}} H(y x^{-1}) f(y) d \mu(y),
	$$
	where $x, y\in \mathbb{G}^n$. If $1 \leq p \leq \infty$  and
	$$
	0 \leq \frac{1}{p} \leq 1,
	$$
	then we have
	$$
	\|\widetilde{T}_H\|_{L^{p} \rightarrow L^p}=\|H\|_{L^1(\mathbb{G}^{n })} .
	$$
\end{cor}
\section{Proof of Theorem \ref{main_2}}

We merely give the proof with the case $m =2$ , and the same is true for the general case $m>2$. Firstly, we prove the sufficiency of Theorem \ref{main_2}. Since $T_2$ is a positive operator, without loss of generality,we take $f  \leq 0$. By the polar coordinates integral formula,we have
\begin{equation}\label{main_11}
	\begin{aligned}
		&\| T_2 (f) \|_{L^p(\mathbb{R}^{n_1}\times \mathbb{R}^{n_2})} \\
		=&\left\{\int_{\mathbb{R}^{n_1}}\int_{\mathbb{R}^{n_2}}\left(\int_{\mathbb{R}^{n_1}}\int_{\mathbb{R}^{n_2}} K(x_1,y_1) K(x_2,y_2) f(y_1,y_2) d y_1 d y_2\right)^p d x_1 d x_2\right\}^{\frac{1}{p}} \\
		=&\left\{\int_{0}^{\infty}\int_{0}^{\infty}\int_{\mathbb{S}^{n_1-1}}\int_{\mathbb{S}^{n_2-1}}\left[ \int_{\mathbb{G}^1}\int_{\mathbb{G}^2}\left( \int_{\mathbb{S}^{n_1-1}}\int{\mathbb{S}^{n_2-1}}K(t_1x_1^{\prime},t_1y_1^{\prime})K(t_2x_2^{\prime},t_2y_2^{\prime})f(t_1y_1^{\prime},t_2y_2^{\prime})\right.\right.\right.\\
		&\left.\left.\left.\times d\sigma(y_1^{\prime})d\sigma(y_2^{\prime}) \right)t_1^{n_1-1}t_2^{n_2-1}d t_1 d t_2\right]^p d\sigma(x_1^{\prime}) d\sigma(x_2^{\prime})t_1^{n_1-1}t_2^{n_2-1} d t_1 d t_2\right\}^{\frac{1}{p}}.
	\end{aligned}
\end{equation}
Let
\begin{equation}\label{12}
	F(f)(t_1, t_2, x_1^{\prime},x_2^{\prime}):=\int_{\mathbb{S}^{n_1-1}} \int_{\mathbb{S}^{n_2-1}}K(t_1 x_1^{\prime},  t_1 y_1^{\prime}) K(t_2 x_2^{\prime},  t_2 y_2^{\prime}) f(t_1 y_1^{\prime},t_2 y_2^{\prime})  d \sigma(y_1^{\prime})d \sigma(y_2^{\prime}).
\end{equation}
Using Minkowski's inequality and (\ref{11}),we have
\begin{equation}\label{15}
	\begin{aligned}
		\| T_2 & (f) \|_{L^p(\mathbb{R}^{n_1}\times \mathbb{R}^{n_2})}\\
		=&\left\{\int_0^{\infty}\int_0^{\infty}\left\{\left[\int_{\mathbb{S}^{n_1-1}}\int_{\mathbb{S}^{n_2-1}}\left(\int_{\mathbb{G}^1}\int_{\mathbb{G}^1} F(f)(t_1,t_2, x_1^{\prime}, x_2^{\prime})  t_1^{n_1-1}t_2^{n_2-1} d t_1 d t_2\right)^p \right.\right.\right.\\
		&\left.\left.\left.\times d \sigma\left(x_1^{\prime}\right)d \sigma(x_2^{\prime})\right]^{\frac{1}{p}}\right\}^p t_1^{n_1-1} t_2^{n_2-1} d t_1 d t_2\right\}^{\frac{1}{p}} \\
		\leq& \left\{\int_0^{\infty}\int_0^{\infty}\left\{\int_{\mathbb{G}^1}\int_{\mathbb{G}^2}\left(\int_{\mathbb{S}^{n_1-1}}\int_{\mathbb{S}^{n_2-1}}\left|F\left(f\right)\left(t_1,t_2, x_1^{\prime},x_2^{\prime}\right)\right|^p d \sigma(x_1^{\prime})d \sigma(x_2^{\prime})\right)^{\frac{1}{p}} \right.\right.\\
		&\left.\left.\times t_1^{n_1-1}t_2^{n_2-1} d t_1 d t_2\right\}^p t_1^{n_1-1}t_2^{n_2-1} d t_1 d t_2\right\}^{\frac{1}{p}}.
	\end{aligned}
\end{equation}
Let $t_1,t_2$ be fixed temporarily.By the definition of $F$ in (\ref{12}), we may regard $F$ as a linear operator from $L^p(\mathbb{S}^{n_1}\times\mathbb{S}^{n_2})$ to $L^p(\mathbb{S}^{n_1}\times\mathbb{S}^{n_2})$.Next , we consider the boundedness of the operator $F$.By the condition of (\ref{11}),$K(t_1 x_1,t_1,y_1) ,  K(t_1,x_1,t_2 y_2)$ as a function of the variables $x_1,x_2,y_1,y_2$ is rotation-invarint.

Define
\begin{equation}\label{main_27}
	A(t_1,t_2)=\int_{\mathbb{S}^{n_1-1}}\int_{\mathbb{S}^{n_2-1}} K(t_1 e_1,t_1 y_1)K(t_2 e_2,t_2 y_2) d\sigma(y_1^{\prime})d\sigma(y_2^{\prime}),
\end{equation}
where $e_i=(1,0),i=1,2$.
By the rotation invariance of $K$, we have
$$
\begin{aligned}
	&\int_{\mathbb{S}^{n_1-1}}\int_{\mathbb{S}^{n_2-1}} K(t_1 e_1,t_1 y_1)K(t_2 x_e,t_2 y_2) d\sigma(y_1^{\prime})d\sigma(y_2^{\prime})\\
	=&\int_{\mathbb{S}^{n_1-1}}\int_{\mathbb{S}^{n_2-1}} K(t_1 x_1,t_1 y_1)K(t_2 x_2,t_2 y_2) d\sigma(y_1^{\prime})d\sigma(y_2^{\prime}).
\end{aligned}
$$
The linear operator $F$ defined in (\ref{12}) satisfies two conditions (1), (2) in Lemma \ref{main_13}. Thus, we have
\begin{equation}\label{16}
	\begin{aligned}
		& \left(\int_{\mathbb{S}^{n_1-1}}\int_{\mathbb{S}^{n_2-1}}\left|F(f)(t_1,  t_2,  x_1^{\prime},x_2^{\prime})\right|^p d \sigma(y_1^{\prime})d \sigma(y_2^{\prime})\right)^{\frac{1}{p}} \\
		\quad \leq & A(t_1,t_2)\left\{\int_{\mathbb{S}^{n_1-1}}\int_{\mathbb{S}^{n_2-1}}\left|f(t_1 y_1^{\prime},t_2 y_2^{\prime})\right|^{p} d \sigma(y_1^{\prime}) d \sigma(y_1^{\prime})\right\}^{\frac{1}{p_i}} .
	\end{aligned}
\end{equation}
Set
$$
h(t_1,t_2) :=\left(\int_{\mathbb{S}^{n_1-1}}\int_{\mathbb{S}^{n_2-1}} |f(t_1y_1,t_2y_2)|^p d\sigma( y_1^{\prime})d \sigma(y_2^{\prime})\right)^\frac{1}{p},
$$
for $0<t_1<\infty,0<t_2<\infty$.
By a computation ,we have
\begin{equation}\label{17}
	\|h\|_{L^{p}(\mathbb{G}^1 \times \mathbb{G}^2)}=\|f\|_{L^{p}(\mathbb{R}^{n_1}\times\mathbb{R}^{n_2})}.
\end{equation}
Combining the inequalities (\ref{15}),(\ref{16}) with (\ref{17}) yields that
\begin{equation}\label{main_21}
	\begin{aligned}
		&\| T_2 (f) \|_{L^p(\mathbb{R}^{n_1}\times\mathbb{R}^{n_2})} \\
		\leq& \left\{\int_0^{\infty}\int_0^{\infty}\left\{\int_{\mathbb{G}^1}\int_{\mathbb{G}^2} A(t_1,t_2) \left\{\int_{\mathbb{S}^{n_1-1}}\int_{\mathbb{S}^{n_2-1}}\left|f(t_1 y_1^{\prime},t_2 y_2^{\prime})\right|^{p} d \sigma(y_1^{\prime}) d \sigma(y_2^{\prime})\right\}^{\frac{1}{p}}\right.\right. \\
		&\left.\left.\times t_1^{n_1-1} t_2^{n_2-1}d t_1d t_2\right\}^p\times t_1^{n_1-1} t_2^{n_2-1}d t_1d t_2\right\}^{\frac{1}{p}} \\
		= & \left\{\int_0^{\infty}\int_0^{\infty}\left[\int_{\mathbb{G}^1}\int_{\mathbb{G}^2}h(t_1,t_2) \int_{\mathbb{S}^{n_1-1}}\int_{\mathbb{S}^{n_2-1}} K(t_1 x_1^{\prime}, t_1 y_1^{\prime}) K(t_2 x_2^{\prime}, t_2y_2^{\prime})d\sigma(y_1^{\prime})d\sigma(y_2^{\prime}) \right.\right.\\ &\left.\left.\times t_1^{\frac{n_1}{p^{\prime}}-1}t_2^{\frac{n_2}{p^{\prime}}-1} d t_1 d t_2\right]^p t_1^{n_1-1}t_2^{n_2-1} d t_1 d t_2 \right\}^{\frac{1}{p}} \\
		= & \left\{\int _ { 0 } ^ { \infty }\int _ { 0 } ^ { \infty } \left[\int_{\mathbb{G}^1}\int_{\mathbb{G}^2}  h(t_1,t_2) \left(\frac{t_1}{t}\right)^{\frac{n_1}{p^{\prime}}}\left(\frac{t_2}{t}\right)^{\frac{n_2}{p^{\prime}}} \int_{\mathbb{S}^{n_1-1}} \int_{\mathbb{S}^{n_2-1}}K\left( x_1^{\prime},\frac{t_1}{t}y_1^{\prime}\right)K\left( x_2^{\prime},\frac{t_2}{t}y_2^{\prime}\right)\right.\right. \\
		& \left.\left.\times  d \sigma\left(y_1^{\prime}\right)d \sigma\left(y_2^{\prime}\right) \frac{d t_1}{t_1}\frac{d t_2}{t_2}\right]^p \frac{d t_1}{t_1}\frac{d t_2}{t_2}\right\}^{\frac{1}{p}},
	\end{aligned}
\end{equation}
where we use the homogeneity of degree $- n$ of the kernel $K$.Notice that the integral
$$
\int_{\mathbb{S}^{n_1-1}} \int_{\mathbb{S}^{n_2-1}}K\left( x_1^{\prime},\frac{t_1}{t}y_1^{\prime}\right)K\left( x_2^{\prime},\frac{t_2}{t}y_2^{\prime}\right)d \sigma(y_1^{\prime})d \sigma(y_2^{\prime})
$$
is independent of $x_1^{\prime},x_2^{\prime}$ because of the rotation of $K$,i.e.,it only depends on $t_1,t_2 $.
Set
\begin{equation}\label{main_18}
	H(t_1, t_2):= t_1^{\frac{n_1}{p^{\prime}}}t_2^{\frac{n_2}{p^{\prime}}} \int_{\mathbb{S}^{n_1-1}} \int_{\mathbb{S}^{n_2-1}}K( x_1^{\prime} ,t_1 y_1^{\prime} ) K( x_2^{\prime}, t_2 y_2^{\prime}) d \sigma(y_1^{\prime})d \sigma(y_1^{\prime}),
\end{equation}
simple computation of integral implies that
$$
\|H\|_{L^1(\mathbb{G}^1\times\mathbb{G}^2)}=\int_{\mathbb{R}^{n_1}} \int_{\mathbb{R}^{n_2}} K(e_1,y_1)K(e_2,y_2) |y_1|^{-\frac{n_1}{p}} |y_2|^{-\frac{n_2}{p}} d y_1 d y_2.
$$
Regarding $H$ in (\ref{main_18}) as a integral kernel,we define a new operator
\begin{equation}\label{main_20}
	\widetilde{T}_H(h_1, h_2)(t):=\int_{\mathbb{G}^1}\int_{\mathbb{G}^2} H\left(\frac{t_1}{t},\frac{t_2}{t}\right) h\left(t_1,t_2\right)\frac{d t_1}{t_1}\frac{d t_2}{t_2},
\end{equation}
where $h(t_1,t_2)$ defined (\ref{17}).We have
\begin{equation}\label{main_22}
	\begin{aligned}
		&\widetilde{T}_H\left(h\right)(t_1,t_2) \\
		\quad=&\int_{\mathbb{G}^1}\int_{\mathbb{G}^2} h\left(t_1,t_2\right) \left(\frac{t_1}{t}\right)^{\frac{n_1}{p^{\prime}}} \left(\frac{t_2}{t}\right)^{\frac{n_2}{p^{\prime}}}\int_{\mathbb{S}^{n_1-1}}\int_{\mathbb{S}^{n_2-1}} K\left( x_1^{\prime},\frac{t_1}{t} y_1^{\prime}\right) K\left( x_2^{\prime},\frac{t_2}{t} y_2^{\prime}\right)\\
		&\times d \sigma\left(y_1^{\prime}\right) d \sigma\left(y_2^{\prime}\right) \frac{d t_1}{t_1}\frac{d t_2}{t_2}.
	\end{aligned}
\end{equation}
It is easy to check that operator $\widetilde{T}_H$ satisfies the all conditions of Corollary \ref{main_19}. Thus it follows from (\ref{17}) and (\ref{main_20}) that
\begin{equation}\label{main_23}
	\begin{aligned}
		\|\widetilde{T}_H(h)\|_{L^p(\mathbb{G}^1\times\mathbb{G}^2)} & \leq\|H\|_{L^1(\mathbb{G}^1\times\mathbb{G}^2)} \|h\|_{L^p{(\mathbb{G}^1\times\mathbb{G}^2)}} \\
		& =\int_{\mathbb{R}^{n_1}}\int_{\mathbb{R}^{n_2}} K(e_1, y_1) K(e_2, y_2) \left|y_1\right|^{-\frac{n_1}{p}} \left|y_2\right|^{-\frac{n_2}{p}}d y_1 d y_2 \left\|f\right\|_{L^{p}\left(\mathbb{R}^{n_1}\times\mathbb{R}^{n_2}\right)} .
	\end{aligned}
\end{equation}
Consequently,it implies from (\ref{main_21}),(\ref{main_22}) and (\ref{main_23}) that
\begin{equation}\label{main_24}
	\begin{aligned}
		\|{T}_2(f)\|_{L^p(\mathbb{R}^{n_1}\times\mathbb{R}^{n_2})} & \leq\|\widetilde{T}_H(h)\|_{L^p(\mathbb{G}^1\times\mathbb{G}^2)} \\
		& \leq\int_{\mathbb{R}^{n_1}}\int_{\mathbb{R}^{n_2}} K(e_1, y_1) K(e_2, y_2) |y_1|^{-\frac{n_1}{p}} |y_2|^{-\frac{n_2}{p}}d y_1 d y_2 \|f\|_{L^{p}(\mathbb{R}^{n_1}\times\mathbb{R}^{n_2})} .
	\end{aligned}
\end{equation}
The inequality (\ref{main_24}) implies that
\begin{equation}\label{main_32}
	\|T_2\|_{L^{p}(\mathbb{R}^{n_1}\times\mathbb{R}^{n_2}) \rightarrow L^p(\mathbb{R}^{n_1}\times\mathbb{R}^{n_2})} \leq \int_{\mathbb{R}^{n_1}} \int_{\mathbb{R}^{n_2}}K(e_1,y_1)K(e_2,y_2) |y_1|^{-\frac{n_1}{p}}|y_2|^{-\frac{n_2}{p}} d y_1 d y_2.
\end{equation}
Next,we will prove the reverse inequality.In fact, it suffices to prove
$$
\|T_2\|_{L^{p}(\mathbb{R}^{n_1}\times\mathbb{R}^{n_2}) \rightarrow L^p(\mathbb{R}^{n_1}\times\mathbb{R}^{n_2})} \geq\|\tilde{T}_H\|_{L^{p}(\mathbb{G}^1\times\mathbb{G}^2) \rightarrow L^p(\mathbb{G}^1\times\mathbb{G}^2)},
$$
where $H$ is defined as(\ref{main_20}).
For any function $h$ : $\mathbb{G}\rightarrow[0,\infty)$, we define the corresponding function $f:\mathbb{R}^n\rightarrow[0,\infty)$ as
$$
f(y_1,y_2)=\omega_{n_1-1}^{-\frac{1}{p}}\omega_{n_2-1}^{-\frac{1}{p}}|y_1|^{-\frac{n_1}{p}} |y_2|^{-\frac{n_2}{p}} h(|y_1|,|y_2|).
$$
Clearly , $f$ is a radial nonnegative function and it easily implies that
\begin{equation}\label{main_26}
	\begin{aligned}
		&\int_{\mathbb{S}^{n_1-1}}\int_{\mathbb{S}^{n_2-1}}K(t_1x_1^{\prime},t_1y_1^{\prime})K(t_2x_2^{\prime},t_2y_2^{\prime})f(t_1 y_1,t_2y_2)d\sigma(y_1)d\sigma(y_2)\\
		=&\omega_{n_1-1}^{-\frac{1}{p}}\omega_{n_2-1}^{-\frac{1}{p}} A(t_1,t_2)t_1^{-\frac{1}{p}}t_2^{-\frac{1}{p}}h(t_1,t_2).
	\end{aligned}
\end{equation}
According to (\ref{main_11}) and (\ref{main_26}) , we have
$$
\begin{aligned}
	& \|T_2(f)\|_{L^p(\mathbb{R}^{n_1}\times\mathbb{R}^{n_2})} \\
	\quad=&\left(\int_0^{\infty}\int_0^{\infty}\int_{\mathbb{S}^{n_1-1}}\int_{\mathbb{S}^{n_2-1}}\left(\int_{\mathbb{G}^1}\int_{\mathbb{G}^2} \omega_{n_1-1}^{-\frac{1}{p}}\omega_{n_2-1}^{-\frac{1}{p}} A(t_1,t_2)  t_1^{-\frac{n_1}{p}} t_2^{-\frac{n_2}{p}} h(t_1,t_2) t_1^{n_1-1} t_2^{n_2-1} d t_1 d t_2\right)^p \right.\\
	&\left.\times d \sigma(x_1^{\prime}) d\sigma(x_2^{\prime}) t_1^{n_1-1} t_2^{n_2-1} d t_1 d t_2\right)^{\frac{1}{p}} \\
	\quad=&\left(\int_0^{\infty}\int_0^{\infty}\left(\int_{\mathbb{G}^1} \int_{\mathbb{G}^2}A(t_1,t_2) h(t_1,t_2) t_1^{\frac{n_1}{p^{\prime}}}t_2^{\frac{n_2}{p^{\prime}}} \frac{d t_1}{t_1}\frac{d t_1}{t_1}\right)^p t_1^{n_1-1} t_2^{n_2-1}d t_1 d t_2\right)^{\frac{1}{p}}.
\end{aligned}
$$
Combining (\ref{main_27}),(\ref{main_18}) with the homogeneity of $-n$, we obtain
$$
A(t_1,t_2)t_1^{\frac{n_1}{p}} t_2^{\frac{n_2}{p}}=H(\frac{t_1}{t},\frac{t_2}{t})t_1^{\frac{n_1}{p}} t_2^{\frac{n_2}{p}}.
$$
Consequently, it follows that
\begin{equation}\label{main_29}
	\|T_2 (f) \|_{L^P (\mathbb{R}^{n_1}\times\mathbb{R}^{n_2})}=\|\widetilde{T}_H(h)\|_{L^P (\mathbb{G}^{1}\times\mathbb{G}^{2})}.
\end{equation}
Since  operator $\widetilde{T}_H$ satisfies all the conditions of Corollary \ref{9}, we have
\begin{equation}\label{28}
	\|\widetilde{T}\|_{L^{p}\rightarrow L^{p}}=\|H\|_{L^1(\mathbb{G}^1\times\mathbb{G}^2)}=\int_{\mathbb{R}^{n_1}} \int_{\mathbb{R}^{n_2}}K(e_1,y_1)K(e_2,y_2) |y_1|^{-\frac{n_1}{p}}|y_2|^{-\frac{n_2}{p}} d y_1 d y_2.
\end{equation}
Since $H$ is nonnegative , we conclude from (\ref{main_29}) that
\begin{equation}\label{29}
	\begin{aligned}
		\|\widetilde{T}_H\|_{L^{p}\rightarrow L^p} =\sup_{\substack{\|h\|}_{L^{p}(\mathbb{G}^1\times\mathbb{G}^2)}\neq0 ,h\geq 0}\frac{\|\widetilde{T}_H (h)\|_{L^p(\mathbb{G}^1\times\mathbb{G}^2)}}{\|h\|_{L^p(\mathbb{G}^1\times\mathbb{G}^2)}}\\
		\\ \leq \sup_{\substack{\|f\|}_{L^{p}(\mathbb{R}^{n_1}\times\mathbb{R}^{n_2})}\neq0}\frac{\|T_2 (f)\|_{L^p(\mathbb{R}^{n_1}\times\mathbb{R}^{n_2})}}{\|f\|_{L^p(\mathbb{R}^{n_1}\times\mathbb{R}^{n_2})}}&
		=\|T_2\|_{L^p(\mathbb{R}^{n_1}\times\mathbb{R}^{n_2})\rightarrow L^p(\mathbb{R}^{n_1}\times\mathbb{R}^{n_2})}.
	\end{aligned}
\end{equation}
It immediately follows from (\ref{28}) and (\ref{29}) that
\begin{equation}\label{main_31}
	\|T_2\|_{L^{p}(\mathbb{R}^{n_1}\times\mathbb{R}^{n_2}) \rightarrow L^p(\mathbb{R}^{n_1}\times\mathbb{R}^{n_2})} \geq\int_{\mathbb{R}^{n_1}} \int_{\mathbb{R}^{n_2}}K(e_1,y_1)K(e_2,y_2) |y_1|^{-\frac{n_1}{p}}|y_2|^{-\frac{n_2}{p}} d y_1 d y_2.
\end{equation}
Both the inequalities (\ref{main_32}) and (\ref{main_31}) , yield
$$
\|T_2\|_{L^{p}(\mathbb{R}^{n_1}\times\mathbb{R}^{n_2}) \rightarrow L^p(\mathbb{R}^{n_1}\times\mathbb{R}^{n_2})} =\int_{\mathbb{R}^{n_1}} \int_{\mathbb{R}^{n_2}}K(e_1,y_1)K(e_2,y_2) |y_1|^{-\frac{n_1}{p}}|y_2|^{-\frac{n_2}{p}} d y_1 d y_2.
$$
So far, we have completed the proof of the Theorem \ref{main_2}.
\section{Sharp constant for  linear Hilbert operator}
In this section, we will consider linear Hilbert operator which can be found in \cite{Gra}. Applying Theorem \ref{main_2}, we can easily obtain the sharp constant on higher-dimensional product spaces.
\begin{defn}
	Let $m \in \mathbb{N}, n_i \in \mathbb{N}, x_i \in \mathbb{R}^{n_i}, 1 \leq i \leq  m$, and f be a nonnegative measurable function on $\mathbb{R}^{n_1}\times\mathbb{R}^{n_2} \times \cdots \times \mathbb{R}^{n_m}$. Then the linear Hilbert operator $T^{*}_m$ on higher-dimensional product spaces is defined by
	\begin{equation}
		T^{*}_mf(x)=\int_{\mathbb{R}^{n_1}}\int_{\mathbb{R}^{n_2}}\cdots\int_{\mathbb{R}^{nm}}\frac{f(y_1,\ldots,y_m)}{\prod_{i=1}^{m}(|x_i|^{n_i}+|y_i|^{n_i})}d y_1 d y_2\cdots d y_m.
	\end{equation}
\end{defn}
\begin{thm}\label{main_111}
	Let $1<p<\infty, m\in\mathbb{N},n_i\in\mathbb{N}, x_i\in\mathbb{R}^{n_i},1\leq i\leq m$. If $f\in L^p(\mathbb{R}^{n_1}\times\mathbb{R}^{n_2}\times\cdots\times\mathbb{R}^{n_m})$, then we have
	$$
	\|T^*_m f\|_{L^p(\mathbb{R}^{n_1}\mathbb{R}^{n_2}\times\cdots\times\mathbb{R}^{n_m})}\leq A.
	$$
	Morever,
	$$
	\|T^*_m f\|_{L^p(\mathbb{R}^{n_1}\times\mathbb{R}^{n_2}\times\cdots\times\mathbb{R}^{n_m})\rightarrow L^p(\mathbb{R}^{n_1}\times\mathbb{R}^{n_2}\times\cdots\times\mathbb{R}^{n_m})}=A,
	$$
	where
	$$
	A=\prod_{i=1}^m\frac{2\pi^{\frac{n_i}{2}}}{\Gamma(\frac{n_i}{2})n_i}m\Gamma(1-\frac{1}{p})\Gamma(\frac{1}{p}).
	$$
\end{thm}
\begin{proof}[Proof of Theorem $\ref{main_111}$]
	Obviously, the operator $T^*_m$ satisfies the condition of Theorem \ref{main_2}, so we have
	$$
	\begin{aligned}
		T^{*}_mf(x)=&\int_{\mathbb{R}^{n_1}}\int_{\mathbb{R}^{n_2}}\cdots\int_{\mathbb{R}^{n_m}}\frac{f(y_1,\ldots,y_m)}{\prod_{i=1}^{m}(|x_i|^{n_i}+|y_i|^{n_i})}dy_1dy_2\cdots d y_m\\
		=&\int_{\mathbb{R}^{n_1}}\int_{\mathbb{R}^{n_2}}\cdots\int_{\mathbb{R}^{nm}}\frac{|y_1|^{-\frac{n_1}{p}}|y_2|^{-\frac{n_1}{p}}\cdots|y_m|^{\frac{n_m}{p}}}{\prod_{i=1}^{m}(1+|y_i|^{n_i})}dy_1dy_2\cdots d y_m.
	\end{aligned}
	$$
	We merely give the proof with the case $m = 2$ , and the same is true for the general case $m > 2$.
	By the polar coordinates integral formula,we have
	$$
	\begin{aligned}
		&\int_{\mathbb{R}^{n_1}}\int_{\mathbb{R}^{n_2}}\frac{|y_1|^{-\frac{n_1}{p}}|y_2|^{-\frac{n_1}{p}}}{(1+|y_1|^{n_1})(1+|y_2|^{n_2})}dy_1d y_2\\
		=&\omega_{n_1-1}\omega_{n_2-1}\int_0^\infty\int^\infty_0\frac{r_1^{\frac{n_1}{p}-1}r_2^{\frac{n_2}{p}-1}}{(1+r_1^{n_1})(1+r_2^{n_2})}d r_1 d r_2\\
		=&\frac{\omega_{n_1-1}}{n_1}\frac{\omega_{n_2-1}}{n_2}\int_0^\infty\int_0^\infty\frac{s_1^{-\frac{1}{p}}s_2^{-\frac{1}{p}}}{(1+s_1)(1+s_2)}d s_1 d s_2.
	\end{aligned}
	$$
	Define
	$$
	F(\beta):=\int_0^\infty\int_0^\infty \frac{s_1^{-\beta}}{1+s_1}\frac{s_2^{-\beta}}{1+s_2} d s_1 d s_2,
	$$
	then we have
	$$
	\begin{aligned}
		F(\beta)&=\int_0^\infty\int_0^\infty\frac{s_1^{-\beta}}{1+s_1}\frac{s_2^{-\beta}}{1+s_2}d s_1 d s_2\\
		&=\int_0^\infty\int_0^\infty\frac{1}{(1+s_1)s_1^{\beta}}\frac{1}{(1+s_2)s_2^{\beta}}d s_1 d s_2\\
		&=\int_0^\infty\frac{1}{(1+s_1)s_1^\beta}d s_1\int_0^\infty\frac{1}{(1+s_2)s_2^\beta}d s_2\\
		&=\int_0^1(1-t_1)^{-\beta}t_1^{\beta-1}d t_1\int_0^1(1-t_2)^{-\beta}t_2^{\beta-1}d t_2\\
		&=2 B(1-\beta,\beta)\\
		&=2\Gamma\left(1-\frac{1}{p}\right)\Gamma\left(\frac{1}{p}\right).
	\end{aligned}
	$$
	Thus, we can obtain
	$$
	\begin{aligned}
		T^*_2f(x)&=\frac{\omega_{n_1-1}}{n_1}\frac{\omega_{n_2-1}}{n_2}2\Gamma\left(1-\frac{1}{p}\right)\Gamma\left(\frac{1}{p}\right)\\
		&=\frac{2\pi^{\frac{n_1}{2}}}{\Gamma\left(\frac{n_1}{2}\right)n_1}\frac{2\pi^{\frac{n_2}{2}}}{\Gamma\left(\frac{n_2}{2}\right)n_2}2\Gamma\left(1-\frac{1}{p}\right)\Gamma\left(\frac{1}{p}\right).
	\end{aligned}
	$$
	This completes the proof of Theorem \ref{main_111}.
\end{proof}
\subsection*{Acknowledgements} This work was supported by National Natural Science Foundation of China (Grant No. 12071473) and Shandong Jianzhu University Foundation (Grant No. XNBS20099).

\begin{flushleft}

	\vspace{0.3cm}\textsc{Xiang Li\\School of Science\\Shandong Jianzhu University\\Jinan, 250000\\P. R. China}
	
	\emph{E-mail address}: \textsf{lixiang162@mails.ucas.ac.cn}
	
		\vspace{0.3cm}\textsc{Zunwei Fu\\Department of Mathmatics\\Linyi University \\Linyi, 270005\\P. R. China}
	
	\emph{E-mail address}: \textsf{fuzunwei@lyu.edu.cn}
	
		\vspace{0.3cm}\textsc{Zhongci Hang\\School of Science\\Shandong Jianzhu University \\Jinan, 250000\\P. R. China}
	
	\emph{E-mail address}: \textsf{babysbreath4fc4@163.com}
	
\end{flushleft}

\end{document}